\documentclass{amsart}

\newcommand{\cal}{\mathcal}
\newcommand{\hslash}{\hbar}

\newcommand{\calg}{$C^*$ -algebra \;}
\newcommand{\calgs}{$C^*$ -algebras \;}

\newcommand{\cpt}{compact\;}

\newcommand{\lcpt}{locally compact\;}
\newcommand{\cnt}{continuous\;}

\newcommand{\fn}{function\;}

\newcommand{\shomo}{$*$-homomorphism\;}

\newcommand{\repn}{representation\;}
\newcommand{\repns}{representations\;}
\newcommand{\pro}{projective limit\;}

\newcommand{\pros}{pro-$C^*$-algebras\;}

\newcommand{\qgs}{quantum groups\;}

\newtheorem{defi}{Definition}[section]
\newtheorem{prop}{Proposition}[section]
\newtheorem{theo}{Theorem}[section]
\newtheorem{lemm}{Lemma}[section]
\newtheorem{cor}{Corollary}[section]

\begin{document}

\title[Generalized Stone's Theorem]{on generalized Stone's theorem}
\author[M. Amini]{Massoud Amini}
\address{Department of Mathematics and Statistics\\ University of Saskatchewan
\\106 Wiggins Road, Saskatoon\\ Saskatchewan, Canada S7N 5E6\\mamini@math.usask.ca}
\thanks{Partially supported by an IPM grant}
\maketitle

\begin{abstract}

It is known that the generator of a strictly continuous one
parameter unitary group in the multiplier algebra of a \calg is
affiliated to that \calg . We show that under natural non degeneracy conditions, this self adjoint
unbounded operator lies indeed in the (unbounded) multiplier algebra
of the Pedersen's ideal of the \calg .
\end{abstract}

\section{Introduction}

One of the main objectives of the Heisenberg formulation of the
Quantum Mechanics is to give appropriate models for the {\it
commutation relations}, the most famous one of which being $$
[P,Q]=i\hslash I $$ where $\hslash$ is the {\it Planck constant}
and $P$, $Q$ are the {\it quantum position} and {\it quantum
momentum}. It was known from the beginning that bounded linear
operators can not satisfy such a relation (convince yourself by
checking this for matrices where you have a {\it trace} for
free!). In particular this does not  happen in a \calg . Although
projective limit of \calgs can include unbounded operators, this
can {\it not} happen in a \pro also. (Just recall that some of
their quotients are \calgs [Ph88b]). One classical trick is to
replace this type of relation by a stronger commutation property,
which in this case is $$ U_t Q U_{-t}=e^{i\hslash t}Q $$ where
$U_t =e^{itH}$, for some closed operator $H$. It has been shown in
[HQV] that any strictly \cnt one parameter group of unitaries
$(U_t )$ in the multiplier algebra of of a \calg $A_0$ has a
generator $H$ ( a version of Stone's theorem) which could be
chosen to be in $A_0 ^{\eta}$ (the set of densely defined linear
operators on $A_0$ which are affiliated with $A_0$ [Wr91]). Here we show
that indeed $H$ could be chosen more specifically. Let
$A_{00}=K(A_0)$ be the Pedersen's ideal of $A_0$ and
$\Gamma(A_{00})$ be the topological algebra of (unbounded)
multiplier on $A_{00}$ [LG]. We show that $H$ could be chosen to
be in $A=\Gamma(A_{00})$.

\section{Stone's theorem}

\bigskip

We start with some lemmas from [Wr91] in which we replace $A_0^{\eta}$
with $A=\Gamma(A_{00})$. We keep the notations of the above paragraph all
over the paper. We also freely use the notations and terminology of [Wr91].
All morphisms are supposed to be non degenerate, and following
[Wr 91] we use $Mor$ to denote the set of morphisms.

The following lemma has been already proved in [Am], but we bring
the proof here for the sake of completeness.

\begin{lemm}
Let $A_0$ and $B_0$ be \calgs . Let $\phi_0 \in Mor(A_0 ,B_0 )$ be
strictly non degenerate (that is $B_{00}\subseteq
\phi_0(A_{00})B$). Then $\phi_0$ extends uniquely to a morphism
$\phi\in Mor(A,B)$ such that $\phi(A_{00})B_0$ is a core for
$\phi(T)$ and $\phi(T)(\phi_0 (a)b)= \phi_0 (Ta)b$ for all $a\in
A_{00}$, $b\in B_0$ , and $ T\in A=\Gamma(A_{00})$. 
\end{lemm}
{\bf Proof}
Since $\phi_0$ preserves the spectral theory,
$\phi_0(A_{00})\subseteq B_{00}$. On the other hand,
$\phi_0(A_{00})$ is clearly a dense ideal of $\phi_0(A_0)$.
Therefore $B_0\phi_0(A_{00})B_0$ is an ideal of $B_0$ which is
dense in $B_0\phi_0(A_0)B_0$. Now if $\phi_0$ is non degenerate
then $\phi_0(A_0)B_0$ is dense in $B_0$, so $B_0\phi_0(A_0)B_0$ is
dense in $B_0 B_0 =B_0$, i.e. $B_0\phi_0(A_{00})B_0$ is a dense
ideal of $B_0$ and so contains $B_{00}$. But
$B_0\phi_0(A_{00})B_0\subseteq B_0 B_{00} B_0\subseteq B_{00}$,
hence the equality holds. If $\phi_0$ is strictly non degenerate
then $\phi_0(A_{00})B_0\supseteq B_{00}$. The converse inclusion
follows from the fact that $\phi_0(A_{00})\subseteq B_{00}$. Hence
$\phi_0(A_{00})B_0=B_{00}$. Now the right hand side is self
adjoint and the adjoint of the left hand side is
$B_0\phi_0(A_{00})$, hence $B_0\phi_0(A_{00})=B_{00}$.\qed

\bigskip

Now each $T\in A=\Gamma(A_{00})$ could be considered
as an element of $A_0^\eta$ and
$\phi_0$ also extends to $\tilde \phi :A_0 ^\eta \to B_0 ^\eta$ by
[Wr91, 1.2]. But there is no ambiguity as we have

\begin{lemm} With the above notations, $\tilde\phi(T)=\phi(T)$. In
particular we have $\phi (z_T )=z_{\phi (T)}$.
\end{lemm}
{\bf Proof} The first statement follows from above lemma and the
facts that by assumption $\phi_0(A_{00})B_0\supseteq B_{00}$ and
$\tilde \phi(T)(\phi_0 (a)b)= \phi_0 (Ta)b\ \ (a\in D(T)\supseteq
A_{00},\ b\in B_0) \ [Wr91,1.2]$. The second statement is proved
for $\tilde\phi$ in [Wr91,1.2], and follows from the first for$\phi$.\qed

\label{mor}
\begin{prop}
Let $A_0$ be a \calg and $T\in A$ be self adjoint. Let $z_T\in M(A_0)$
be the $z$-transform of $T$. Then $\sigma(z)\subseteq [-1,1]$. Assume
that
$$
span\{f(z_T)a: f\in C_{00}(-1,1), a\in A_0\}\supseteq A_{00}.
$$
Then there is a
unique $\phi=\phi_T \in Mor(C(\Bbb R ), A)$ such that $\phi_T
(id)=T$. Moreover $\phi (C_{00}(\Bbb R ))A_0\supseteq A_{00}$ and
$\phi_T (z)=z_T$, for $z(t)=t(1+t^2)^{\frac{1}{2}}$,\  $t\in \Bbb R$.
\end{prop}
{\bf Proof} $z_T\in M(A_0 )$ is self adjoint and $\|z_T\|\leq 1$,
so its spectrum $\sigma(z_T )$ is contained in $[-1,1]$. The same
is true for $z\in C_b (\Bbb R )$. We use the continuous functional
calculus to show the uniqueness. If $\phi (id)=T$ then
$\phi(z)=z_T$ and $\phi(f\circ z)=f(z_T)$, for each $f\in C_b
(-1,1)=C[-1,1]$. But each element of $C_0 (\Bbb R )$ is of the
form $f\circ z$, where $f\in C_0 (-1,1)$ and the uniqueness
follows.

For the existence, let's define $\phi_0 :C_0(\Bbb R)\to M(A_0)$ by
$\phi_0 (f\circ z)=f(z_T)$. Then $\phi_0 $ is clearly a \shomo . Also
by assumption, $\phi_0 (C_{00}(\Bbb R)A_0)\supseteq A_{00}$
which means that $\phi_0$ is strictly non degenerate.
By [Am, Theorem 3.4],
$\phi_0$ extends (uniquely) to some $\phi\in Mor(C(\Bbb
R),\Gamma(A_{00}))$. Now $\phi(f\circ z)=f(z_T)$, for each $f\in
C_0(-1,1)$. In particular, for $f(t)=t(1-t^2)^{\frac{-1}{2}}$, we get
$f\circ z=id$ and so $\phi(id)=f(z_T)=T$.
\qed

Let $\phi_0$ and $\phi$ be as in the proof of the above proposition. Then
$\phi$ is not injective in general. But we can do the classical trick
to make it injective. Consider
$ker(\phi_0)=\{f\in C_0(\Bbb R): f=0 \ \text{on}\ \sigma(T)\}$. Then
$C_0(\Bbb R)/ker(\phi_0)=C_0(\sigma(T))$ and the corresponding
quotient map identifies with the restriction map $\pi :C_0(\Bbb R)\to
C_0(\sigma(T))$. Let $id_{\sigma(T)}=\pi(id)$, then we have

\begin{prop}
Let $A_0$ be a \calg and $T\in \Gamma(A_{00})$, then there is a unique
embedding $\psi=\psi_T\in Mor(C(\sigma(T)),\Gamma(A_{00}))$ such that
$\psi_T (id_{\sigma(T)})=T$.
\end{prop}
{\bf Proof} Since $\pi$ is onto, there is a function $\psi_0$ such
that $\psi_0\pi=\phi_0$. Then $\psi_0$ is a \shomo . To see that
it is strictly non degenerate it is enough to observe that $\psi_0
(C_0(\sigma(T)))A_0 =\phi_0(C_0(\Bbb R))A_0\supseteq A_{00}$.
Therefore it extends to a morphism $\psi\in
Mor(C(\sigma(T)),\Gamma (A_{00}))$.

Next let us show that $\psi$ is one-one. It is clear that $\psi_0$ is
one-one (sine $\psi_0 =\phi_0\pi$ where $\pi(ker\phi_0 )=\{0\}$). Take
any $F\in C(\Bbb R)$ such that $\psi(F)=0$, then
$0=\psi(F)\psi_0(f)a=\psi_0(Ff)a$, for each $f\in C_0(\sigma(T))$ and
$a\in A_0$. This means that $\psi_0(Ff)\in M(A_0)$ multiplies $A_0$
into $0$, i.e. $\psi_0(Ff)=0$. Hence $Ff=0$ for each $f\in
C_0(\sigma(T))$, and so $F=0$.
\qed

\bigskip

Now we are prepared to prove the generalization of Stone's
theorem. Let $A_0$ be a \calg , $A_{00}=K(A_0)$ its Pedersen's
ideal, and $A=\Gamma(A_{00})$ be the algebra of (unbounded)
multipliers of $A_{00}$ [LT]. For each $t\in \Bbb R$ consider the
function $e_t \in C(\Bbb R)$ defined by $e_t (s)=exp(its)\ \ (s\in
\Bbb R)$. Let $h\in A=\Gamma(A_{00})$ and $U_t
=\phi_h(e_t)=exp(ith)$. Then $(U_t)_{t\in \Bbb R}$ is a one
parameter strictly \cnt unitary group in $M(A_0)$ (the strict
continuity follows from the fact that $\phi_h$ is non degenerate).
Moreover, if $h\in b(\Gamma(A_{00}))=M(A_0)$, then this is also
norm continuous. Conversely each strictly \cnt unitary group in
$M(A_0)$ is of this form, for some $h\eta A_0$ [HQV]. Here we want
$h$ to be actually in $\Gamma(A_{00})$. Clearly for this to
happen, we would need to put some condition on the unitary group.
This is the content of the following result. The proof is quite
similar to [HQV, 2.1]. Here we only sketch those parts of the
proof which have to be modified. But first a definition.

\begin{defi}
Let $A_0$ be a \calg and $(U_t)_{t\in \Bbb R}$ be a one parameter
strictly \cnt unitary group in $M(A_0)$. Let's define
$\alpha :L^1 (\Bbb R)\to M(A_0)$ by $\alpha(f)=\int_{\Bbb R}
f(t)U_t dt$, where
$$
\bigg(\int f(t)U_t dt\bigg)x=\int f(t)U_t x dt\ \ (x\in A_0)
$$
is in the Bochner sense. This extends to a \shomo $\alpha\in Mor(C^*
(\Bbb R), A_0)$ [HQV]. We say that $(U_t)$ is (strictly) non degenerate,
if the morphism $\alpha$ is (strictly) non degenerate
(cf. [Am]).
\end{defi}

\label{stone}
\begin{theo}{\bf (Generalized Stone's Theorem)}
Let $A_0$ be a \calg and $(U_t)_{t\in \Bbb R}$ be a one parameter
strictly \cnt strictly non degenerate unitary group in $M(A_0)$.
Then there exists a self
adjoint $h\in \Gamma(K(A_0)))$ such that $U_t =exp(ith)$ for $t\in
\Bbb R$. Moreover, if $(U_t)_{t\in \Bbb R}$ is norm continuous, then
$h\in M(A_0)$.
\end{theo}
{\bf Proof} Define $\alpha\in Mor(C^* (\Bbb R), A_0)$ as above.
Then $\alpha$ is strictly non degenerate, and so it extends to a
morphism of the corresponding \pros , which we still denote it by
$\alpha\in Mor(\Gamma(K(C^* (\Bbb R))), \Gamma(K(A_0)))$ (see the
discussion before Theorem 3.2 in [Am]). Also it is well known that
the Fourier transform $\cal{F} :L^1 (\Bbb R)\to C_0(\Bbb R)$
extends to an isomorphism $\cal F \in Mor(C^* (\Bbb R), C_0(\Bbb
R))$ with $\cal F (e_t)=\lambda_t$, where $\lambda$ is the left
regular \repn of $\Bbb R$. Since $\cal F$ is surjective we can
extend it to an isomorphism $\cal F\in Mor(\Gamma(K(C^* (\Bbb
R))), C(\Bbb R))$ [LT]. Also it is clear that
$\alpha(\lambda_t)=U_t\ (t\in \Bbb R)$. Define $h=(\alpha\circ
\cal F ^{-1})(id)\in \Gamma (A_{00})$. This is self adjoint and by
the uniqueness part of Proposition 1.2 we have $\phi_h
=\alpha\circ \cal F ^{-1}$. Now $U_t
=\alpha(\lambda_t)=\alpha(\cal F
^{-1}(e_t))=\phi_h(e_t)=exp(ith)$, which finishes the proof of the
first part of the theorem.

Now, for each $t\in \Bbb R$ we get
\begin{align*}
\|U_t -1\|
 &=\|\phi_h(e_t -1)\|=\|\psi_h\pi(e_t -1)\|=\|\pi(e_t -1)\|\\
 &=sup\{|e^{it\lambda} -1|:\lambda\in \sigma (h)\}.
\end{align*}
So, if the unitary group is norm \cnt , then $\sigma(h)$ is
bounded. Hence $h\in b(\Gamma(K(A_0)))$ [Ph88b]. But
$b(\Gamma(K(A_0)))=M(A_0)$ [Ph88a], so $h\in M(A_0)$ and the
proof is complete.
\qed

\bigskip

Next, following [HQV],we show that $h$ can be found by {\it
differentiating} the unitary group $(U_t)_{t\in \Bbb R}$.

\begin{prop}
Let $A_0$ be a \calg ,$(U_t)_{t\in \Bbb R}$ a one parameter
strictly \cnt unitary group in $M(A_0)$, and $h$ a self
adjoint element of $ \Gamma(K(A_0)))$ such that $U_t =exp(ith)$ for $t\in
\Bbb R$. Define the (unbounded) operator $H:D(H)\subseteq A_0\to A_0$ by
\begin{gather}
D(H)=\{a\in A_0 : t\mapsto U_t a \ \text{is} \ C^1 \} \\
Ha=\frac{d}{dt}|_{t=0} U_t a=\lim_{t\to 0} (U_t a-a)/t \ \ (a\in D(H)).
\end{gather}
Then $h\subseteq iH=i\frac{d}{dt}|_{t=0} U_t$.
\end{prop}
{\bf Proof} Using minimality of $A_{00}$ we have $A_{00}\subseteq
\sqrt{1-z_h} A_0$. Now this last set is contained in $D(H)=D(iH)$
[HQV, 2.2]. The fact that $h(x)=iH(x)$, for each $x\in A_{00}$,
follows from the calculations of [HQV, 2.2].\qed

\bigskip

This in particular shows the uniqueness of the element $h\in
\Gamma(A_{00})$ for which $U_t=exp(ith)\ \ (t\in\Bbb R)$. We call this
element the {\it infinitesimal generator} of the one parameter group
$(U_t)_{t\in\Bbb R}$.

\section{elements affiliated with group \calg}

The group \calgs are important objects in the theory of quantum
groups. One reason is that they are {\it dual objects} to continuous
functions. When the underlying group is not discrete (non compact
case), one would expect some unbounded elements to come into the
play. These can not belong to the group $C^*$-algebra, but they are usually
affiliated with it. The problem of finding all elements affiliated
with $C^{*}(G)$ is open in general. Some attempts are done to find them
in the case that $G$ is a Lie group [Wr91]. In this case we know
that the elements of the corresponding Lie algebra,
considered as differential operators, are affiliated with the group
\calg. The question of whether all elements of the universal
enveloping algebra of $G$ are affiliated with $C^{*}(G)$ was left
open. In this
section we give a partial answer to this question using
\pros.

Let $G$ be a Lie group, and $\pi:G\to
B(\cal{X})$ be a (strongly \cnt) \repn of $G$ on a Banach space
$\cal{X}$. We say that $x\in\cal{X}$ is a $C^{\infty}$-{\it
vector}({\it analytic vector}, respectively) of $\pi$ if the map
$g\mapsto\pi(g)x$ from $G$ to $\cal{X}$ is a $C^{\infty}$
(analytic, respectively) \fn. We denote the set of all such elements
$x\in\cal{X}$ by
$D^{\infty}=D^{\infty}(\pi)$ ($D^{\omega}(\pi)$, respectively). Then
this is a dense subset of $\cal{X}$. Indeed L. G\aa rding showed
that if $\phi\in C_{00}^{\infty}(G)$ and
$\pi(\phi)x=\int_{G} \pi(t)x\phi(t)dt \ \ (x\in\cal{X})$, then
$D_{00}^{\infty}=\pi(C_{00}^{\infty}(G))\cal{X}\subseteq\cal{X}$
is dense. This is called the {\it G\aa rding domain} of $\pi$. Let
$\cal G=L(G)$ be the Lie algebra of $G$, i.e. the set of all
left invariant vector fields on $G$ (at identity). For each
$X\in\cal G$, define
$$
\pi(X)x=\lim_{h\to 0} \frac{\pi(\exp{hX})x-x}{h}\ \ \ (x\in D^{\infty})
$$
Then
$\pi(X):D^{\infty}\subseteq\cal{X}\to D^{\infty}\subseteq\cal{X}$
is a densely defined unbounded operator on $\cal{X}$, and it
 is {\it skew symmetric} if $\pi$ is unitary (and $\cal{X}$ a Hilbert
 space).

Harish-Chandra noticed that $D^{\infty}$ could have a
subspace $D$ such that $\pi(X)D \subseteq D$ but
$\pi(g)\bar{D}\not\subseteq\bar{D}$. Therefore he suggested replacing
$D^{\infty}$ with $D^{\omega}$.  He showed that $D^{\omega}$ is dense
in $\cal X$ for certain  \repns of a semisimple Lie group [Har]. P.
Cartier and J. Dixmier  gave a proof  for all unitary \repns [CD2], and
E. Nelson used a generalization  of the
fundamental solution of the heat
equation on Lie groups  to prove this for an arbitrary \repn.
He also used analytic  vectors to give a sufficient condition for a
\repn of a Lie  algebra to be induced by a unitary \repn of the Lie
group [Nel]. Later  I. Segal showed that this is equivalent to the {\it
complete  positivity} of the \repn with respect to an appropriate
cone  [Seg].

We are mainly interested in the case of the universal
\repn.   For any \lcpt group $G$, the universal \repn $u: G\to
M(C^{*}(G))$  is determined by the following universal property:
 For each \calg $A$ and each \repn $\pi:G\to A$, there is a unique
$\tau\in Mor(C^{*}(G), A)$ such that $\pi=\tau u$. Then $u$ is  continuous
and open with respect to the strict topology.   Now let $G$ be a Lie
group and let $\cal G$ be its Lie algebra of  dimension $N$ with a basis
$X_1,X_2,\dots,X_N$.  Then the {\it universal
enveloping algebra} $\cal U=\cal U (\cal G)$  of $\cal G$ is an
$*$-algebra under $X^{*}=-X$, $X\in\cal G$.  Elements of $\cal U$
are {\it differential operators} on $G$  commuting with the  right
translations: Take $D ^{\infty}= D^{\infty}(u)$, then each  $X\in\cal
U$ defines
\begin{gather*}
 du(X): D^{\infty}\subseteq C^{*}(G)\to B(H_{u})  \\
 du(X)a=Xu(g)a|_{g=e}=\lim_{t\to 0} \frac{1}{t}(u(\exp{tX})a- a)
\end{gather*}
which is a closable operator whose closure is simply denoted  by $X:
D^{\infty}\subseteq C^{*}(G)\to C^{*}(G)$.   For the reasons explained
in the beginning of this section,  we are interested in elements
affiliated with $C^{*}(G)$. We  know that for $X\in\cal U$, if the
differential equation $X^ {*}Xf=-f$ has only trivial bounded
$C^{\infty}$-solution, then  $X,X^{*}\eta C^{*}(G)$ [Wr95,2.1]. It is
known that all the elements of $\cal G$ and also the {\it
elliptic operator}  $\Delta=-\sum_{1}^{N} X_{i}^{2}$, where $X_i$'s
form a basis of $\cal G$, are affiliated
with $C^{*}(G)$ [Wr95].  The
$X_i$'s are skew self-adjoint, where as $\Delta$ is
self-adjoint  and positive (yes positive!).

To look more carefully at the above problem, one has to distinguish
between the right and left multipliers. To each $g\in G$ we associate
a bounded operator $u(g)\in B(H_u)$. We can view this as a bounded
multiplier on $C^* (G)$, by considering it as an element of $C^*
(G)^{**}$. This means that we put $u(g)=(S_g ,T_g )$, where $S_g$ and
$T_g$ are bounded extensions of the mappings
$$
S_g(u(f))=u(\delta_g *f)\ \ \text{and}\ \ T_g(u(f))=u(f*\delta_g)\ \
(f\in L^1 (G))
$$
Now consider the following
$$
D_\ell ^{\infty}=\{a\in C^* (G): g\mapsto S_g (a) \ is \ a \ C^{\infty} \
map\}
$$
and define
\begin{gather*}
 d_\ell u(X): D_\ell ^{\infty}\subseteq C^{*}(G)\to  C^{*}(G)  \\
 d_\ell u(X)a=\lim_{t\to 0} \frac{1}{t}(S_{\exp{(tX)}}(a)- a)\ \ \
(a\in  D_\ell ^{\infty})
\end{gather*}
then $D_\ell ^{\infty}$ is a dense right ideal of  $C^{*}(G)$. In order
to connect this with the above problem, we need to define the one sided
version of some of the previous notions.

\begin{defi} Let $J$ be a right ideal in a \calg $A_0$. Then a left
multiplier of $J$ is a linear map $S:J\to J$ such that $S(ax)=S(a)x\ \
(a\in J, x\in A_0 )$. The set of all such maps is denoted by $\Gamma_\ell
(J)$. A right multiplier is defined similarly and we denote the set of
all right multipliers by $\Gamma_r (J)$.
\end{defi}

There is a natural left version of the $\kappa$-topology for $\Gamma_\ell
(J)$: A net $(S_\alpha )\in \Gamma_\ell (J)$ converges to $S\in \Gamma_\ell
(J)$ iff $\| S_\alpha(a)-S(a)\|\to 0$. We call this the
$\ell\kappa$-topology. Obviously $J$ is $\ell\kappa$-dense in $\Gamma_\ell
(J)$. Similarly on can define $r \kappa$-topology on $\Gamma_r (J)$.

\begin{lemm} If $J\subseteq A_0$ is a dense right ideal, then $b(\Gamma_\ell
(J))\subseteq M_\ell (A_0 )$.
\end{lemm}{\bf Proof} This is to say that each bounded left multiplier of
$J$ extends uniquely to a left multiplier of $A_0$, which follows
from density of $J$ in $A_0$. \qed

\begin{defi}
We say that $T\in L(A_0 )$ is left-affiliated with $A_0$ if $z_T=T(I+T^*
T)^{\frac{-1}{2}}\in M_\ell (A_0 )$ and $z_T ^* z_T < I$, or
equivalently, if there is $z\in M_\ell (A_0 )$ of norm one such that
$\sqrt{1-z_T ^* z_T} A_0$ is dense in $A_0$ and $T( \sqrt{1-z_T ^*
z_T} a)=za$. We denote the set of all these elements by
$A_0 ^{\ell\eta}$. The elements right-affiliated with $A_0$ are defined
similarly, and the set of all such elements is denoted by
$A^{r\eta}$.
\end{defi}

The following lemma is now trivial.

\begin{lemm} For each dense right ideal $J\in A_0$, $\Gamma_\eta (J)\subseteq
A_0 ^{\ell\eta}$.
\qed
\end{lemm}

Now let us come back to the special case of $A_0 =C^* (G)$ and $J=D_\ell
^{\infty}$. Take any $X\in \cal U (\cal G)$, then for each $a\in J$,
$ d_\ell u(X)a=\lim_{t\to 0} \frac{1}{t}(S_{\exp{(tX)}}-1)(a)$, which
means that $d_\ell u(X)= \ell\kappa$-$lim \ S_{\exp{(tX)}}-1 \in
\Gamma_\ell (J)\subseteq
A_0 ^{\ell\eta}$. Therefore we have

\begin{theo}
Let $G$ be a Lie group and $\cal G$ be its Lie algebra. For each
$X\in \cal U (\cal G)$, $d_r u(X)$ and $d_\ell u(X)$ are
respectively right-affiliated and left affiliated with $C^* (G)$.
\qed
\end{theo}

We can get a better result using the generalized Stone's theorem,
proved in previous section. Let $A_0$ be a \calg and $G$ be a Lie
group acting on $A_0$ through a strictly \cnt \repn $u:G\to
M(A_0)$. Then
$$
C^{\infty}(u)=\{x\in A_0 : g\mapsto u(g)x \ \text{is}\ C^{\infty}\}
$$
contains the dense subspace of $A_0$ spanned by all elements $u(f)x$,
with $f\in C_{00}^{\infty}(G)$ and $x\in A_0$.

\begin{prop}
With the above notation, to each $X\in \cal G$ there corresponds an
element $\hat X\in\Gamma(A_{00})$ such that
$$
e^{t\hat X}=u(e^{tX})\ \ (t\in\Bbb R ).
$$
Moreover $\hat X$ leaves $C^{\infty}(u)$ invariant and
$[X,Y]\ \hat{}=[\hat X,\hat Y]$, when restricted to  $C^{\infty}(u)$,
for all $X,Y\in \cal G$.
\end{prop}{\bf Proof} Take $X\in \cal G$ and put $U_t =u(e^{tX})$. Then
Theorem 2.1 applies and gives $\hat X\in \Gamma(A_{00})$ with $U_t
=e^{t\hat X}$ for each $t\in\Bbb R  $ (here $\hat X$ is {\it skew}
self adjoint). The rest is proved as in [HQV, 2.4]. \qed

\label{Ualg}
\begin{cor}
With the above notation, $\cal U (\cal G )\subseteq \Gamma(K(C^*
(G)))$.
\end{cor}
{\bf Proof} Just observe that $G$ acts strictly continuously on
$C^* (G)$ through the universal \repn and apply above proposition.
\qed

\bigskip

{\small {\bf Acknowledgement}: This paper is part of the author's
Ph.D. thesis in the University of Illinois at Urbana-Champaign
under the supervision of Professor Zhong-Jin Ruan. I would like to
thank him for his moral support and scientific guidance during my
studies.}

\end{document}